\documentclass
{Uzhgorod-Mathematical-Paper-2021}
\umpnumber{41} \umpsubnumber{2} \umpyear{2022} \umpissnprint{2616-7700}
\umpissnonline{2708-9568}
 \umpdoi{10.24144/2616-7700}

\usepackage{mathtools,dsfont}

\begin{document}

\selectlanguage{english}
\udk{519.21}
\author{M.~Ilienko}
\institution{National Technical University of Ukraine \lq\lq Igor Sikorsky Kyiv Polytechnic Institute\rq\rq, Kyiv}
\position{associate prof. of the department of mathematical analysis and probability theory}
\academictitle{phD, associate prof.}
\mail{mari-run@ukr.net}
\orcid{0000-0001-6982-6124}

\author{A.~Polishchuk}
\institution{Institute of Applied Physics, National Academy of Science of Ukrine, Kyiv}
\position{junior researcher}
\mail{nastya.varennikova312@gmail.com}

\shortauthor{Ilienko~M., Polishchuk~A.}
\capitalauthor{M.~ILIENKO, A. POLISHCHUK}
\authoreng{Ilienko~M.~K., Polishchuk~A.Yu.}

\title{On the convergence of Baum-Katz series for sums of linear 2-nd order autoregressive sequences of random variables}
\shorttitle{BAUM-KATZ SERIES FOR 2-ND ORDER AUTOREGRESSIVE SEQUENCES \ldots}
\titleeng{On the convergence of Baum-Katz series for sums of linear 2-nd order autoregressive sequences}
\maketitle
\received{.}

\newcommand{\hm}[1]{#1\nobreak\discretionary{}{\hbox{\ensuremath{#1}}}{}}
\newcommand{\Leb}{\mathsf{Leb}}

\begin{abstract}
We consider complete convergence and closely related Hsu-Robbins-Erd\H{o}s-Spitzer-Baum-Katz series for sums whose terms are elements of a linear 2-nd order autoregressive sequences of random variables and prove sufficient conditions for the convergence of this series.

	
\end{abstract}

\noindent
\paragraph{Introduction}
On a common probability space $(\Omega, \mathcal{F},\mathbb{P})$ consider a linear 2-nd order autoregressive sequence of random variables (r.v.'s) $(\xi_{k}, \, k\geq1)$, which obeys the system of stochastic recurrence equations
\begin{equation}\label{model}
\xi_{-1}=0, \ \ \xi_{0}=0, \ \  \xi_{k}=a\xi_{k-1}+b\xi_{k-2}+\theta_{k}, \ \ k\ge 1,
\end{equation}
where $a$ and $b$ are some real constants, $(\theta_{k})$ is a sequence of independent copies of a r.v. $\theta$.
Note that linear regression models of different types have been studied for a long time. A multitude of publications contain various problems for regression sequences of r.v.'s and their extensions. See, for instance \cite{b,mikosh}, and numerous references therein.

For elements of the sequence \eqref{model} set
\begin{equation*}
S_n=\sum_{k=1}^n \xi_k, \ \ n\geq1,
\end{equation*}
and for any $\varepsilon>0$ consider the following series
\begin{equation}\label{BKseries}
\sum_{n=1}^\infty n^{\frac{r}{p}-2}\mathbb{P} \Big\{\frac{|S_n|}{n^{1\!/\!p}}>\varepsilon  \Big\},
\end{equation}
where  $0<p<2$ and $r\geq p$. In this paper we are interested in conditions for the convergence of this series.
Hereinafter we will refer to the series \eqref{BKseries} as to Baum-Katz series, although some other no less prominent authors were involved in introducing it.

Historically, for the sequence $(X_{n}, \, n\geq1)$ of independent copies of a r.v. $X$, and $S_n=\sum_{k=1}^n X_k,$ $n\geq1$, the reduced version (the case $r=2p=2$) of the series \eqref{BKseries} initially arose in the paper by Hsu and Robbins along with the notion of complete convergence, see \cite{Hsu}. In their paper, authors proved that if $\mathbb E X^{2}<\infty$, then  $\sum_{n=1}^\infty\mathbb{P} \Big\{\Big|S_n/n-\mathbb EX\Big|>\varepsilon  \Big\} <\infty$, while the converse was provided by Erd\H{o}s, see \cite{Erd}.
Note that in view of the Borel-Cantelli Lemma complete convergence implies almost sure convergence, and so is tightly connected with the Strong Law of Large Numbers.
Further, Spitzer, see \cite{Sp}, showed that
$\sum_{n=1}^\infty n^{-1}\mathbb{P} \Big\{\Big|S_n/n-\mathbb EX\Big|>\varepsilon  \Big\} <\infty$ if and only if $\mathbb E|X|<\infty$. Note that series \eqref{BKseries} covers the Spitzer's case with $r=p=1$. Finally, for the sequence of independent copies of a r.v. $X$, Baum and Katz, see \cite{BK}, introduced the series \eqref{BKseries} and proved that it is convergent if and only if $\mathbb E|X|^{r}<\infty,$ with $\mathbb E X=0$, when  $r\geq1.$ Since then these classical results have been generalized in several directions, including Banach space setting (see, e.g., \cite{Huan19}). We refer to \cite{deli} where more detailed history on the topic is provided.
Among all extensions we distinguish results concerning complete convergence and convergence of Baum-Katz series for weighted sums of independent r.v.'s, also known as rowwise independent random arrays, (see, e.g., \cite{deli,HuMoritz,gut,Hu,Cai,ChenMaSung} and references therein).

As to dependent patterns, in the paper \cite{ilienko21} necessary and sufficient conditions for the convergence of the series (\ref{BKseries}) were
obtained for the case of first-order autoregressive sequence of r.v.'s, i.e. with $b=0$ in \eqref{model}.

Specifically, in this paper we concentrate on sufficient conditions of the series \eqref{BKseries} for sums of elements of model \eqref{model}, and under some simple assumptions imposed on $a$ and $b$ we expect to obtain similar to independent case Baum-Katz result. In our investigation we intend to reduce the case to the idea provided in \cite{ilienko21}, which in its turn was partially borrowed from \cite{gut}.

\paragraph{Preliminaries}
Consider the nonrandom recurrence sequence $(u_{n} , n\geq 1)$:
\begin{equation}\label{uuu}
u_{-1}=0, \ \ u_{0}=1, \ \ u_{n}=au_{n-1}+bu_{n-2}, \ \ n\geq 1,
\end{equation}

Evaluating \eqref{model} one has
\begin{equation*}
\xi_1=\theta_1=u_0\theta_1,
\end{equation*}
\begin{equation*}
\xi_2=a\theta_1+\theta_2=u_1\theta_1+u_0\theta_2,
\end{equation*}
\begin{equation*}
\xi_3=(a^2+b)\theta_1+a\theta_2+\theta_3=u_2\theta_1+u_1\theta_2+u_0\theta_3,
\end{equation*}
\begin{equation*}
\ldots ,
\end{equation*}
i.e.
\begin{equation*}
\xi_k=\sum_{l=1}^k u_{k-l}\theta_l, \  \ k\geq 1.
\end{equation*}

Now, for $n\geq 1 $ and $1\leq k \leq n$ set
\begin{equation*}
u(n-k)=\sum_{m=0}^{n-k} u_m.
\end{equation*}

Thus,
\begin{equation*}
S_n=\sum_{k=1}^n \xi_k = \sum_{k=1}^n \big( \sum_{l=1}^k u_{k-l}\theta_l \Big)=
\sum_{k=1}^n \Big(\sum_{m=0}^{n-k} u_m \Big) \theta_k=\sum_{k=1}^n u(n-k)\theta_k, \  \ n\geq1.
\end{equation*}

\paragraph{Main result}
Let us immediately proceed to the main result of this paper.

\begin{theorem}\label{theorem1}
	Let in \eqref{model}
\begin{equation}\label{coef_assump}
-1<b<1-|a|,
\end{equation}
and $0<p<2$, $r\geq p$. If  $\mathbb E|\theta|^{r}<\infty$, where $\mathbb E \theta=0 $ whenever $r\geq1$, then for any $\varepsilon>0$,
	\begin{equation*}
	\sum_{n=1}^\infty n^{\frac{r}{p}-2}\mathbb{P} \Big\{\frac{|S_n|}{n^{1\!/\!p}}>\varepsilon  \Big\} <\infty.
	\end{equation*}
		
\end{theorem}

\paragraph{Proof}
In \cite{ilienko21} the analogue of Theorem 1 for linear 1-st order autoregressive sequence of r.v.'s was proved in all details.
We may adopt the proof of sufficiency to our case if we show that the values $u(n-k)$, $n\geq 1 $, $1\leq k \leq n$, are bounded in a similar way (compared with $a(n,k)$'s in \cite{ilienko21}).

Introduce two real-valued matrices
$$
M=\begin{pmatrix}
1 & 0  \\
0 & 0  \\
\end{pmatrix},
\quad
C=\begin{pmatrix}
a & b \\
1 & 0 \\
\end{pmatrix}.
$$
Note, that $C$ is a Frobenius matrix. Let $\lambda_1$ and $\lambda_2$ be its eigenvalues, i.e. roots of the characteristic equation $\lambda^2-a\lambda-b=0$. Denote by $\nu_1$ and $\nu_2$ multiplicities of $\lambda_1$ and $\lambda_2$ respectively.
Set
$$
\rho=\max\big\{|\lambda_1|,|\lambda_2|\big\} \quad \textrm{and} \quad
\mu= \underset{1\leq k\leq 2}{\max} \{\nu_k: |\lambda_k|=\rho\}.
$$
Obviously, in our case either $\mu=1$ or $\mu=2$. Moreover, assumption \eqref{coef_assump} implies that both roots $\lambda_1$ and $\lambda_2$ lie within the unit circle, that is $\rho<1$.

Observe that
\begin{equation*}
CM=
\begin{pmatrix}
a & 0 \\
1 & 0 \\
\end{pmatrix}=\begin{pmatrix}
u_1 & 0 \\
u_0 & 0 \\
\end{pmatrix},
\ \
C^2 M=
\begin{pmatrix}
a^2+b & 0 \\
a & 0 \\
\end{pmatrix}=\begin{pmatrix}
u_2 & 0 \\
u_1 & 0 \\
\end{pmatrix},
\end{equation*}
Further,
\begin{equation*}
C^3 M=
\begin{pmatrix}
a^3+2ab & 0 \\
a^2+b & 0 \\
\end{pmatrix}=\begin{pmatrix}
u_3 & 0 \\
u_2 & 0 \\
\end{pmatrix},
\end{equation*}
and so on.
Using the method of mathematical induction it is easy to show that for any $s\geq 1$,
\begin{equation*}
C^s M=
\begin{pmatrix}
u_s & 0 \\
u_{s-1} & 0 \\
\end{pmatrix}.
\end{equation*}

Let for a square matrix $A=(a_{ij})_{i,j=1}^2$ with real entries $\|\cdot\|$ denote the matrix norm of the following form:
$\|A\|=\Big(\sum_{i,j=1}^2 a_{i,j}^2\Big)^{1/2}.$
According to result by Koval', see \cite{ko} (see also Lemma 7.7.3 \cite{b}), if $C$ is a Frobenius matrix, there exist some constants $c_2>c_1>0$, such that for any $s\geq1$,
\begin{equation*}
c_1 \cdot \rho^s \cdot s^{\mu-1} \leq \|C^s M\|\leq c_2 \cdot \rho^s \cdot s^{\mu-1},
\end{equation*}
where $c_1$ and $c_2$ do not depend on $s$.

Now since
\begin{equation*}
(C^0+C^1+C^2+\ldots+C^{n-k})M=
\begin{pmatrix}
u(n-k) & 0 \\
u(n-k-1)  & 0 \\
\end{pmatrix},
\end{equation*}
then according to Koval's result,
\begin{align}
	 |u(n-k)|\leq & \sqrt{(u(n-k))^2+(u(n-k-1))^2} = \Big\|\sum_{l=0}^{n-k} C^{l}M \Big\| \leq \sum_{l=0}^{n-k}  \Big\|C^{l}M \Big\| \leq \nonumber
	\\
	& \leq c_2\sum_{l=0}^{n-k}  \rho^l \cdot l^{\mu-1},
\label{inequal}
	\end{align}
where $c_2$ is some positive constant.

Now distinguish between two cases:

1) $\lambda_1\neq\lambda_2$ (if $b\neq -a^2/4$). In this case $\mu=1$ and according to \eqref{inequal},
\begin{equation*}
|u(n-k)|\leq c_2\sum_{l=0}^{n-k}  \rho^l= c_2 \frac{1-\rho^{n-k+1}}{1-\rho},  \ \  1\leq k \leq n,  \ \ n\geq 1.
\end{equation*}

2) $\lambda_1=\lambda_2$ (if $b=-a^2/4$). In this case $\mu=2$ and according to \eqref{inequal},
\begin{equation*}
|u(n-k)|\leq c_2\sum_{l=1}^{n-k}  l \rho^l\leq c_2\sum_{l=1}^{\infty}  l \rho^l = c_2\frac{\rho}{(1-\rho)^2},  \ \  1\leq k \leq n,  \ \ n\geq 1.
\end{equation*}

Combining both cases,
\begin{equation*}
|u(n-k)|\leq L=const= \max\Big\{ \frac{c_2}{1-\rho}, \frac{c_2\rho}{(1-\rho)^2} \Big\}, \ \  1\leq k \leq n,  \ \ n\geq 1.
\end{equation*}

Now, briefly adopt the proof of the sufficiency of Theorem 1 in \cite{ilienko21} to our case.
As in \cite{ilienko21}  we first restrict our proof to the case of symmetrically distributed r.v. $\theta$.
	Let us fix any $\varepsilon>0$ and apply an iteration of the Hoffmann-J$\o$rgensen inequality (see \cite{gut} or \cite{ilienko21}) with $s=t=n^{1/p}\varepsilon$. Thus, for $j\geq1$ there exist some constants $C_j$ and $D_j$ such that
	\begin{align}
	&\mathbb{P}\Big\{|S_n|>n^{1/p}\varepsilon\cdot 3^j\Big\} \leq \nonumber
	\\
	&C_j \sum_{k=1}^n \mathbb{P}\Big\{ \Big|u(n-k)\theta_k \Big|>n^{1/p}\varepsilon\Big\}+D_j\Big(\mathbb{P}\Big\{|S_n|>n^{1/p}\varepsilon\Big\} \Big)^{2^j}.  \label{HJ1}
	\end{align}
	
The first terms in (\ref{HJ1}) can be estimated as follows
	\begin{align*}
	&\sum_{k=1}^n \mathbb{P} \Big\{\Big|u(n-k) \theta_k\Big|>n^{1/p}\varepsilon\Big\}=
	\sum_{k=1}^n \mathbb{P} \Big\{|\theta_k|>\frac{n^{1/p}\varepsilon}{|u(n-k)|} \Big\}\leq
	\\
	&\leq\sum_{k=1}^n \mathbb{P} \Big\{|\theta_k|>n^{1/p}\varepsilon L^{-1}\Big\}=n\mathbb{P} \Big\{|\theta|>n^{1/p}\varepsilon_2\Big\},
	\end{align*}
	where $\varepsilon_2= \varepsilon L^{-1}$. Further, we refer to the corresponding estimations in \cite{ilienko21}.

Now consider the second term in (\ref{HJ1}). According to Markov inequality for $r>p$, one has
	\begin{equation*}
	\mathbb{P}\Big\{|S_n|>n^{1/p}\varepsilon\Big\}\leq \frac{\mathbb E|S_n|^r}{(n^{1/p}\varepsilon)^r}.
	\end{equation*}
	Next deal with $\mathbb E|S_n|^r$ distinguishing between the following cases.
	
	\textbf{1)} Let $0<r \leq 1$. Applying the $c_r$-inequality (see, for example, \cite{LBai}) with $c_r=1$ to $\mathbb E|S_n|^r$, one obtains
	\begin{equation*}
	\mathbb E|S_n|^r \leq
	\sum_{k=1}^{n} \mathbb E\Big|u(n-k) \theta_k \Big|^r=
	\sum_{k=1}^{n} \big|u(n-k) \big|^r \mathbb E |\theta_k|^r\leq
	\mathbb E|\theta|^r L^r n.
	\end{equation*}	
	
\textbf{2)} Let $r>1$. In this case to $\mathbb E|S_n|^r$ we consequently apply the Marcinkiewicz-Zygmund inequality (see, for example, \cite{LBai}) and the following well-known inequality: for positive $a_i$, $1\leq i \leq n$,  $n \in \mathbb N$ and $r>0$ it is true that
	\begin{equation*}
	(a_1^2+a_2^2+...+a_n^2)^{r/2}\leq n^{0\vee (r/2-1)} \sum_{i=1}^n a_i^{r}.
	\end{equation*}
	Thus,
	\begin{align*}
	&\mathbb E|S_n|^r\leq
	b_r \mathbb E\Big(\sum_{k=1}^{n} \Big(u(n-k) \theta_k \Big)^2 \Big)^{r/2} \leq
	b_r  n^{0\vee (r/2-1)}  \mathbb E \sum_{k=1}^{n} \big|u(n-k) \theta_k \big|^{r} =
	\\
	&=b_r  n^{0\vee (r/2-1)}  \mathbb E|\theta|^{r}  \sum_{k=1}^{n} \big|u(n-k)  \big|^{r} \leq
	b_r  n^{0\vee (r/2-1)}  \mathbb E|\theta|^{r} L^r n=b_r  n^{1\vee (r/2)}  \mathbb E|\theta|^{r} L^r.
	\end{align*}
Here $b_r$ is some positive constant from the Marcinkiewicz-Zygmund inequality.

Combining the above two cases, we arrive at the following bounds
	\begin{equation*}
	\mathbb E|S_n|^r \leq
	C(r) \mathbb E|\theta|^r {n^{1 \vee (r/2)}},
	\end{equation*}
	with $C(r)=L^r$ or $b_r L^r$ depending on whether $0<r \leq 1$ or $r>1$.

Now to finish the proof one needs to literary follow the steps of it in \cite{ilienko21}.

\begin{example}
If $\theta$ is a normally distributed r.v. with $\mathbb E \theta=0$, the model \eqref{model} represents the so-called Gaussian 2-Markov sequence of r.v.'s. with constant coefficients. In this case the series \eqref{BKseries} converges provided that $0<p<2$, $r\geq p$ and $-1<b<1-|a|$.
\end{example}

\paragraph{Conclusions}
In the paper for sequences of sums whose terms are elements of 2-nd order linear autoregressive sequences, sufficient conditions for the convergence of Baum--Katz series are considered. Under some anticipated assumptions imposed on the coefficients of autoregressive sequence, obtained sufficient conditions are expressed as moment assumption of the generating r.v. The latter, in its turn, agrees with the classical Baum--Katz independent case.

We intently focused our attention on 2-nd order autoregressive sequences, evading general $m$-th order case, since for $m=2$ assumptions imposed on the coefficients of the sequence are described in the most simple form. But, in prospect, by means of the same technique set problem may be generalized to $m$-th order autoregressive sequences for any $m\geq 2$. Moreover, we expect to prove also necessary conditions for convergence of Baum--Katz series for such sequences.

\selectlanguage{ukrainian}
\umpmakeengtitle

\selectlanguage{ukrainian}

\end{document}